\documentstyle[12pt]{article}

\font\twlgot =eufm10 scaled \magstep1 \font\egtgot =eufm8
\font\sevgot =eufm7 \font\twlmsb =msbm10 scaled \magstep1
\font\egtmsb =msbm8 \font\sevmsb =msbm7

\newfam\gotfam

\textfont\gotfam\twlgot \scriptfont\gotfam\egtgot
\scriptscriptfont\gotfam\sevgot

\newfam\msbfam
\textfont\msbfam\twlmsb \scriptfont\msbfam\egtmsb
\scriptscriptfont\msbfam\sevmsb
\def\Bbb{\protect\pBbb}
\def\pBbb{\relax\ifmmode\expandafter\Bb\else\typeout{You cann't use
Bbb in text mode}\fi}
\def\Bb #1{{\fam\msbfam\relax#1}}

\def\thebibliography#1{\bigskip\section*{}\bigskip\list
{$^{\arabic{enumi}}$}{\settowidth\labelwidth{#1}\leftmargin\labelwidth
\advance\leftmargin\labelsep
\usecounter{enumi}}
\def\newblock{\hskip .11em plus .33em minus .07em}
\sloppy\clubpenalty4000\widowpenalty4000 \sfcode`\.=1000\relax}

\def\op#1{\mathop{\fam0 #1}\limits}

\newcommand{\beq}{\begin{equation}}
\newcommand{\eeq}{\end{equation}}
\newcommand{\ben}{\begin{eqnarray}}
\newcommand{\een}{\end{eqnarray}}
\newcommand{\be}{\begin{eqnarray*}}
\newcommand{\ee}{\end{eqnarray*}}
\newcommand{\bea}{\begin{eqalph}}
\newcommand{\eea}{\end{eqalph}}

\newcommand{\cH}{{\cal H}}
\newcommand{\cF}{{\cal F}}

\newcommand{\cG}{{\cal G}}

\newcommand{\bs}{{\bf s}}

\newcommand{\al}{\alpha}

\newcommand{\bt}{\beta}
\newcommand{\dl}{\delta}
\newcommand{\la}{\lambda}

\newcommand{\Om}{\Omega}
\newcommand{\m}{\mu}

\newcommand{\g}{\gamma}

\newcommand{\vt}{\vartheta}
\newcommand{\vf}{\varphi}

\newcommand{\di}{{\rm dim\,}}

\newcommand{\si}{\sigma}

\newcommand{\w}{\wedge}

\newcommand{\dr}{\partial}
\newcommand{\ar}{\op\longrightarrow}
\newcommand{\ot}{\otimes}

\newcounter{eqalph}
\newcounter{equationa}
\newcounter{remark}
\newcounter{example}
\newcounter{theorem}
\newcounter{proposition}
\newcounter{lemma}
\newcounter{corollary}
\newcounter{definition}
\def\theremark{\arabic{remark}}
\def\thetheorem{\arabic{theorem}}

\def\thedefinition{\arabic{definition}}

\newenvironment{proof}{\noindent
{\it Proof:} }{ \medskip}

\newenvironment{theo}{\refstepcounter{theorem}
\bigskip\noindent{\it Theorem \thetheorem:} }{\medskip}

\newenvironment{lem}{\refstepcounter{theorem}
\bigskip\noindent{\it Lemma \thetheorem:}}{\medskip}

\newenvironment{defi}{\refstepcounter{definition}
\bigskip\noindent{\it Definition \thedefinition:}}{\medskip}

\newenvironment{eqalph}{\stepcounter{equation}
\setcounter{equationa}{\value{equation}} \setcounter{equation}{0}

\begin{eqnarray}}{\end{eqnarray}
\setcounter{equation}{\value{equationa}}}

\newcommand{\mar}[1]{}

\begin{document}
\hbox{}

{\parindent=0pt

{\large \bf Global action-angle coordinates for completely
integrable systems with noncompact invariant submanifolds}
\bigskip

{\sc E.Fiorani}\footnote{Electronic mail:
emanuele.fiorani@unicam.it}

{\sl Department of Mathematics and Informatics, University of
Camerino, 62032 Camerino (MC), Italy}

\medskip

{\sc G.Sardanashvily}\footnote{Electronic mail:
sard@grav.phys.msu.su}

{\sl Department of Theoretical Physics, Moscow State University,
117234 Moscow, Russia}

\bigskip
\bigskip

The obstruction to the existence of global action-angle
coordinates of Abelian and noncommutative (non-Abelian) completely
integrable systems with compact invariant submanifolds has been
studied. We extend this analysis to the case of noncompact
invariant submanifolds.

 }

\bigskip
\bigskip

\noindent {\bf I. INTRODUCTION}
\bigskip

We consider Abelian and noncommutative (non-Abelian) completely
integrable Hamiltonian systems (henceforth CISs) on symplectic
manifolds. The Liouville--Arnold (or Liouville--Mineur--Arnold)
theorem for Abelian CISs$^{1-4}$  and the Mishchenko--Fomenko
theorem for noncommutative ones$^{5-8}$ state the existence of
action-angle coordinates around a compact invariant submanifold of
a CIS. These theorems have been extended to the case of noncompact
invariant submanifolds.$^{9-12}$ In particular, this is the case
of time-dependent CISs.$^{13-14}$ Any time-dependent CIS of $m$
degrees of freedom can be represented as the autonomous one of
$m+1$ degrees of freedom on a homogeneous momentum phase space,
where time is a generalized angle coordinate. Therefore, we
further consider only autonomous CISs.

If invariant submanifolds of a CIS are compact, a topological
obstruction to the existence of global action-angle coordinates
has been analyzed.$^{8,15,16}$ Here, we aim extending this
analysis to the case of noncompact invariant submanifolds
(Theorems \ref{nc6}, \ref{nc6'} and \ref{nc0}).

Throughout the paper, all functions and maps are smooth, and
symplectic manifolds are real smooth and paracompact. We are not
concerned with the real-analytic case because a paracompact
real-analytic manifold admits the partition of unity by smooth,
not analytic functions. As a consequence, sheaves of modules over
real-analytic functions need not be acyclic that is essential for
our consideration.

\begin{defi} \label{i0} \mar{i0}
Let $(Z,\Om)$ be a $2n$-dimensional connected symplectic manifold,
and let $(C^\infty(Z), \{,\})$ be the Poisson algebra of smooth
real functions on $Z$. A subset $H=(H_1,\ldots,H_k)$, $n\leq
k<2n$, of $C^\infty(Z)$ is called a (noncommutative) CIS if the
following conditions hold.

(i) All the functions $H_i$ are independent, i.e., the $k$-form
$\op\w^kdH_i$ nowhere vanishes on $Z$. It follows that the map
$H:Z\to \Bbb R^k$ is a submersion, i.e.,
\mar{nc4}\beq
H:Z\to N=H(Z) \label{nc4}
\eeq
is a fibered manifold over a connected open subset $N\subset\Bbb
R^k$.

(ii) There exist smooth real functions $s_{ij}$ on $N$ such that
\mar{nc1}\beq
\{H_i,H_j\}= s_{ij}\circ H, \qquad i,j=1,\ldots, k. \label{nc1}
\eeq

(iii) The matrix function $\bs$ with the entries $s_{ij}$
(\ref{nc1}) is of constant corank $m=2n-k$ at all points of $N$.
\end{defi}

In Hamiltonian mechanics, one can think of the functions $H_i$ as
being integrals of motion of a CIS which are in involution with
its Hamiltonian. Their level surfaces (fibers of $H$) are
invariant submanifolds of a CIS.

If $k=n$, then $\bs=0$, and we are in the case of an Abelian CIS.
If $k>n$, the matrix $\bs$ is necessarily non-zero, and a CIS is
said to be noncommutative.

Note that, in many physical models, the condition (i) of
Definition \ref{i0} fails to hold. In a general setting, one
supposes that the subset $Z_R\subset Z$ of regular points, where
$\op\w^kdH_i\neq 0$, is open and dense. Then one considers a CIS
on this subset. However, a CIS on $Z_R$ fail to be equivalent to
the original one because there is no morphism of Poisson algebras
$C^\infty(Z_R)\to C^\infty(Z)$. In particular, canonical
quantization of the Poisson algebra $C^\infty(Z_R)$, e.g., with
respect to action-angle variables essentially differs from that of
$C^\infty(Z)$.$^{17-19}$ For instance, let $M$ be a connected
compact invariant manifold of an Abelian CIS through a regular
point $z\in Z_R\subset Z$. There exists its open saturated
neighborhood $U_M\subset Z_R$ (i.e., a fiber of $H$ through a
point of $U_M$ belongs to $U_M$) which is a trivial fiber bundle
in tori. By virtue of the above mentioned Liouville--Arnold
theorem, $U_M$ is provided with the Darboux action-angle
coordinates. Then one treats quantization of the Poisson algebra
$C^\infty(U_M)$ with respect to these coordinates as quantization
'around' an invariant submanifold $M$.

Given a CIS in accordance with Definition \ref{i0}, the above
mentioned generalization of the Mishchenko--Fomenko theorem to
noncompact invariant submanifolds states the following.$^{12}$

\begin{theo} \label{nc0'} \mar{nc0'}
Let the Hamiltonian vector fields $\vt_i$ of the functions $H_i$
be complete, and let the fibers of the fibered manifold $H$
(\ref{nc4}) be connected and mutually diffeomorphic. Then the
following hold.

(I) The fibers of $H$ (\ref{nc4}) are diffeomorphic to a toroidal
cylinder
\mar{g120}\beq
\Bbb R^{m-r}\times T^r. \label{g120}
\eeq

(II) Given a fiber $M$ of $H$ (\ref{nc4}), there exists an open
saturated neighborhood $U_M$ of it which is a trivial principal
bundle with the structure group (\ref{g120}).

(III) The neighborhood $U_M$ is provided with the bundle
(generalized action-angle) coordinates $(I_\la,p_A,q^A, y'^\la)$,
$\la=1,\ldots, m$, $A=1,\ldots,n-m$, where $(y'^\la)$ are
coordinates on a toroidal cylinder, such that the symplectic form
$\Om$ on $U_M$ reads
\be
\Om= dI_\la\w dy'^\la + dp_A\w dq^A,
\ee
and a Hamiltonian of a CIS is a smooth function only of the action
coordinates $I_\la$.
\end{theo}

Theorem \ref{nc0'} restarts the Mishchenko--Fomenko one if its
condition is replaced with that the fibers of the fibered manifold
$H$ (\ref{nc4}) are compact and connected.

The proof of Theorem \ref{nc0'} is based on the following
facts.$^{8,12}$ Any function constant on fibers of the fibration
$H$ (\ref{nc4}) is the pull-back of some function on its base $N$.
Due to item (ii) of Definition \ref{i0}, the Poisson bracket
$\{f,f'\}$ of any two functions $f,f'\in C^\infty(Z)$ constant on
fibers of $H$ is also of this type. Consequently, the base $N$ of
$H$ is provided with a unique coinduced Poisson structure
$\{,\}_N$ such that $H$ is a Poisson morphism.$^{20}$ By virtue of
condition (iii) of Definition \ref{i0}, the rank of this coinduced
Poisson structure equals $2(n-m)=2\di N-\di Z$. Furthermore, one
can show the following.$^{8,21}$

\begin{lem} \label{nc8} \mar{nc8} The fibers of the
fibration $H$ (\ref{nc4}) are maximal integral manifolds of the
involutive distribution spanned by the Hamiltonian vector fields
of the pull-back $H^*C$ of Casimir functions $C$ of the coinduced
Poisson structure on $N$.
\end{lem}

In particular, a Hamiltonian of a CIS is the pull-back onto $Z$ of
some Casimir function of the coinduced Poisson structure on $N$.

It follows from Lemma \ref{nc8} that invariant submanifolds of a
noncommutative CIS are maximal integral manifolds of a certain
Abelian partially integrable system (henceforth PIS).

\begin{defi} \label{i1} \mar{i1} A collection $\{S_1,\ldots, S_m\}$
of $m\leq n$ independent smooth real functions in involution on a
$2n$-dimensional symplectic manifold $(Z,\Om)$ is called a PIS.
\end{defi}

Let us consider the map
\mar{g106}\beq
S:Z\to W\subset\Bbb R^m. \label{g106}
\eeq
Since functions $S_\la$ are everywhere independent, this map is a
submersion onto an open subset $W\subset \Bbb R^m$, i.e., $S$
(\ref{g106}) is a fibered manifold of fiber dimension $2n-m$.
Hamiltonian vector fields $v_\la$ of functions $S_\la$ are
mutually commutative and independent. Consequently, they span an
$m$-dimensional involutive distribution on $Z$ whose maximal
integral manifolds constitute a foliation $\cF$ of $Z$. Because
functions $S_\la$ are constant on leaves of this foliation, each
fiber of a fibered manifold $Z\to W$ (\ref{g106}) is foliated by
the leaves of the foliation $\cF$. If $m=n$, we are in the case of
an Abelian CIS, and the leaves of $\cF$ are connected components
of fibers of the fibered manifold (\ref{g106}). The
Poincar\'e--Lyapounov--Nekhoroshev theorem$^{22-24}$ generalizes
the Liouville--Arnold one to a PIS if leaves of the foliation
$\cF$ are compact. It imposes a sufficient condition which
Hamiltonian vector fields $v_\la$ must satisfy in order that the
foliation $\cF$ is a fibered manifold.$^{24,25}$ Extending the
Poincar\'e--Lyapounov--Nekhoroshev theorem to the case of
noncompact integral submanifolds, we in fact assumed from the
beginning that these submanifolds formed a fibration.$^{11,14,18}$
Here, we aim to prove the following global variant of Theorem 6 in
Ref. [11].

\begin{theo} \label{nc6} \mar{nc6}
Let a PIS $\{S_1,\ldots,S_m\}$ on a symplectic manifold $(Z,\Om)$
satisfy the following conditions.

(i) The Hamiltonian vector fields $v_\la$ of $S_\la$ are complete.

(ii) The foliation $\cF$ is a fiber bundle $\cF:Z\to N$.

(iii) Its base $N$ is simply connected.

(iv) The cohomology $H^2(N,\Bbb Z)$ of $N$ with coefficients in
the constant sheaf $\Bbb Z$ is trivial.

\noindent Then the following hold.

(I) The fiber bundle $\cF$ is a trivial principal bundle with the
structure group (\ref{g120}), and we have a composite fibered
manifold
\mar{g107}\beq
S=\zeta\circ\cF: Z\ar N\ar W, \label{g107}
\eeq
where $N\to W$ however need not be a fiber bundle.

(II) The fibered manifold (\ref{g107}) is provided with the
adapted fibered (generalized action-angle) coordinates
\be
(I_\la,x^A,y'^\la)\to (I_\la,x^A)\to (I_\la), \qquad
\la=1,\ldots,m, \quad A=1,\ldots 2(n-m),
\ee
such that the coordinates $(I_\la)$ possess identity transition
functions, and the symplectic form $\Om$ reads
\mar{nc3'}\beq
\Om= dI_\la\w dy'^\la + \Om_A^\la dI_\la\w dx^A+ \Om_{AB} dx^A\w
dx^B. \label{nc3'}
\eeq
\end{theo}

If one supposes from the beginning that leaves of the foliation
$\cF$ are compact, the condition (i) of Theorem \ref{nc6} always
holds, and the assumption (ii) can be replaced with the
requirement that $\cF$ is a fibered manifold with mutually
diffeomorphic connected fibers. Recall that any fibered manifold
whose fibers are diffeomorphic either to $\Bbb R^r$ or a compact
connected manifold $K$ is a fiber bundle.$^{26}$ However, a
fibered manifold whose fibers are diffeomorphic to a product $\Bbb
R^r\times K$ (e.g., a toroidal cylinder) need not be a fiber
bundle (see Ref. [27], Example 1.2.2).

Theorem \ref{nc6} is proved in Section II. Since $m$-dimensional
fibers of the fiber bundle $\cF$ admit $m$ complete independent
vector fields, they are locally affine manifolds diffeomorphic to
a toroidal cylinder (\ref{g120}). Then the condition (iii) of
Theorem \ref{nc6} guaranties that the fiber bundle $\cF$ is a
principal bundle with the structure group (\ref{g120}).
Furthermore, this principal bundle is trivial due to the condition
(iv), and it is provided with the bundle action-angle coordinates.
Note that conditions (ii) and (iii) of Theorem \ref{nc6} are
sufficient, but not necessary.

 If $m=n$, the following corollary of Theorem \ref{nc6} states the
existence of global action-angle coordinates of an Abelian CIS.

\begin{theo} \label{nc6'} \mar{nc6'}
Let an Abelian CIS $\{H_1,\ldots, H_n\}$ on a symplectic manifold
$(Z,\Om)$ satisfy the following conditions.

(i) The Hamiltonian vector fields $\vt_i$ of $H_i$ are complete.

(ii) The fibered manifold $H$ (\ref{nc4}) is a fiber bundle with
connected fibers over a simply connected base $N$ whose cohomology
$H^2(N,\Bbb Z)$ is trivial.

\noindent Then the following hold.

(I) The fiber bundle $H$ (\ref{nc4}) is a trivial principal bundle
with the structure group (\ref{g120}).

(II) The symplectic manifold $Z$ is provided with the global
Darboux coordinates $(I_\la,y'^\la)$ such that $\Om= dI_\la\w
dy'^\la$.
\end{theo}

Due to Lemma \ref{nc8}, a manifested global generalization of
Theorem \ref{nc0'} is a corollary of Theorem \ref{nc6} (see
Section III).

\begin{theo} \label{nc0} \mar{nc0}
Given a noncommutative CIS in accordance with Definition \ref{i0},
let us assume the following.

(i) Hamiltonian vector fields $\vt_i$ of integrals of motion $H_i$
are complete.

(ii) The fibration $H$ (\ref{nc4}) is a fiber bundle with
connected fibers.

(iii) Let $V$ be an open subset of the base $N$ of this fiber
bundle which admits $m$ independent Casimir functions of the
coinduced Poisson structure on $N$.

(iv) Let $V$ be simply connected, and let the cohomology
$H^2(V,\Bbb Z)$ be trivial.

\noindent Then the following hold.

(I) The fibers of $H$ (\ref{nc4}) are diffeomorphic to a toroidal
cylinder (\ref{g120}).

(II) The restriction $Z_V$ of the fiber bundle $H$ (\ref{nc4}) to
$V$ is a trivial principal bundle with the structure group
(\ref{g120}).

(III) The fiber bundle $Z_V$ is provided with the bundle
(generalized action-angle) coordinates $(I_\la,x^A,y'^\la)$ such
that the action-angle coordinates $(I_\la,y'^\la)$ possess
identity transition functions and the symplectic form $\Om$ on
$Z_V$ reads
\mar{nc3}\beq
\Om= dI_\la\w dy'^\la + \Om_{AB} dx^A\w dx^B. \label{nc3}
\eeq
\end{theo}

Note that, if invariant submanifolds of a CIS are assumed to be
connected and compact, condition (i) of Theorem \ref{nc0} is
unnecessary since vector fields $v_\la$ on compact fibers of $H$
are complete. In this case, condition (ii) of Theorem \ref{nc0}
also holds because, as was mentioned above, a fibred manifold with
compact mutually diffeomorphic fibers is a fiber bundle.

In the case of an Abelian CIS, the coinduced Poisson structure on
$N$ equals zero, the integrals of motion $H_\la$ are the pull-back
of $n$ independent functions on $N$, and Theorem \ref{nc0} reduces
to Theorem \ref{nc6'}.

Following the original Mishchenko--Fomenko theorem, let us mention
noncommutative CISs whose integrals of motion $\{H_1,\ldots,H_k\}$
form a $k$-dimensional real Lie algebra $\cG$ of rank $m$ with the
commutation relations
\be
\{H_i,H_j\}= c_{ij}^h H_h, \qquad c_{ij}^h={\rm const.}
\ee
In this case, complete Hamiltonian vector fields $\vt_i$ of $H_i$
define a locally free Hamiltonian action on $Z$ of some simply
connected Lie group $G$ whose Lie algebra is isomorphic to
$\cG$.$^{28,29}$ Orbits of $G$ coincide with $k$-dimensional
maximal integral manifolds of the regular distribution on $Z$
spanned by Hamiltonian vector fields $\vt_i$.$^{30}$ Furthermore,
one can treat $H$ (\ref{nc4}) as an equivariant momentum mapping
of $Z$ to the Lie coalgebra $\cG^*$, provided with the coordinates
$x_i(H(z))=H_i(z)$, $z\in Z$.$^{18,31}$ In this case, the
coinduced Poisson structure $\{,\}_N$ coincides with the canonical
Lie--Poisson structure on $\cG^*$ given by the Poisson bivector
field
\be
w=\frac12 c_{ij}^h x_h\dr^i\w\dr^j.
\ee
Casimir functions of the Lie--Poisson structure are exactly the
coadjoint invariant functions on $\cG^*$. They are constant on
orbits of the coadjoint action of $G$ on $\cG^*$ which coincide
with leaves of the symplectic foliation of $\cG^*$. Let $V$ be an
open subset of $\cG^*$ which obeys the conditions (iii) and (iv)
of Theorem \ref{nc0}. Then the open subset $H^{-1}(V)\subset Z$ is
provided with the action-angle coordinates.

\bigskip
\bigskip

\noindent {\bf II. PROOF OF THEOREM 3}
\bigskip

In accordance with the well-known theorem,$^{28,29}$ complete
Hamiltonian vector fields $v_\la$ define an action of a simply
connected Lie group on $Z$. Because vector fields $v_\la$ are
mutually commutative, it is the additive group $\Bbb R^m$ whose
group space is coordinated by parameters $s^\la$ with respect to
the basis $\{e_\la=v_\la\}$ for its Lie algebra. The orbits of the
group $\Bbb R^m$ in $Z$ coincide with the fibers of the fiber
bundle
\mar{g105}\beq
\cF:Z\to N. \label{g105}
\eeq
Since vector fields $v_\la$ are independent on $Z$, the action of
$\Bbb R^m$ on $Z$ is locally free, i.e., isotropy groups of points
of $Z$ are discrete subgroups of the group $\Bbb R^m$. Given a
point $x\in N$, the action of $\Bbb R^m$ on a fiber
$M_x=\cF^{-1}(x)$ factorizes as
\mar{d4}\beq
\Bbb R^m\times M_x\to G_x\times M_x\to M_x \label{d4}
\eeq
through the free transitive action of the factor group $G_x=\Bbb
R^m/K_x$, where $K_x$ is the isotropy group of an arbitrary point
of $M_x$. It is the same group for all points of $M_x$ because
$\Bbb R^m$ is a commutative group. Since the fibers $M_x$ are
mutually diffeomorphic, all isotropy groups $K_x$ are isomorphic
to the group $\Bbb Z_r$ for some fixed $0\leq r\leq m$.
Accordingly, the groups $G_x$ are isomorphic to the Abelian group
\mar{g100} \beq
G=\Bbb R^{m-r}\times T^r, \label{g100}
\eeq
and fibers of the fiber bundle (\ref{g105}) are diffeomorphic to
the toroidal cylinder (\ref{g100}).

Let us bring the fiber bundle (\ref{g105}) into a principal bundle
with the structure group (\ref{g100}). Generators of each isotropy
subgroup $K_x$ of $\Bbb R^m$ are given by $r$ linearly independent
vectors $u_i(x)$ of the group space $\Bbb R^m$. These vectors are
assembled into an $r$-fold covering $K\to N$. This is a subbundle
of the trivial bundle
\mar{g101}\beq
N\times R^m\to N \label{g101}
\eeq
whose local sections are local smooth sections of the fiber bundle
(\ref{g101}). Such a section over an open neighborhood of a point
$x\in N$ is given by a unique local solution $s^\la(x')e_\la$ of
the equation
\be
g(s^\la)\si(x')=\exp(s^\la v_\la)\si(x')=\si(x'), \qquad
s^\la(x)e_\la=u_i(x),
\ee
where $\si$ is an arbitrary local section of the fiber bundle
$Z\to N$ over an open neighborhood of $x$. Since $N$ is simply
connected, the covering $K\to N$ admits $r$ everywhere different
global sections $u_i$ which are global smooth sections
$u_i(x)=u^\la_i(x)e_\la$ of the fiber bundle (\ref{g101}). Let us
fix a point of $N$ further denoted by $\{0\}$. One can determine
linear combinations of the functions $S_\la$, say again $S_\la$,
such that $u_i(0)=e_i$, $i=m-r,\ldots,m$, and the group $G_0$ is
identified to the group $G$ (\ref{g100}). Let $E_x$ denote the
$r$-dimensional subspace of $\Bbb R^m$ passing through the points
$u_1(x),\ldots,u_r(x)$. The spaces $E_x$, $x\in N$, constitute an
$r$-dimensional subbundle $E\to N$ of the trivial bundle
(\ref{g101}). Moreover, the latter is split into the Whitney sum
of vector bundles $E\oplus E'$, where $E'_x=\Bbb R^m/E_x$.$^{32}$
Then there is a global smooth section $\g$ of the trivial
principal bundle $N\times GL(m,\Bbb R)\to N$ such that $\g(x)$ is
a morphism of $E_0$ onto $E_x$, where $u_i(x)=\g(x)(e_i)=\g_i^\la
e_\la$. This morphism is also an automorphism of the group $\Bbb
R^m$ sending $K_0$ onto $K_x$. Therefore, it provides a group
isomorphism $\rho_x: G_0\to G_x$. With these isomorphisms, one can
define the fiberwise action of the group $G_0$ on $Z$ given by the
law
\mar{d5}\beq
G_0\times M_x\to\rho_x(G_0)\times M_x\to M_x. \label{d5}
\eeq
Namely, let an element of the group $G_0$ be the coset
$g(s^\la)/K_0$ of an element $g(s^\la)$ of the group $\Bbb R^m$.
Then it acts on $M_x$ by the rule (\ref{d5}) just as the coset
$g((\g(x)^{-1})^\la_\bt s^\bt)/K_x$ of an element
$g((\g(x)^{-1})^\la_\bt s^\bt)$ of $\Bbb R^m$ does. Since entries
of the matrix $\g$ are smooth functions on $N$, the action
(\ref{d5}) of the group $G_0$ on $Z$ is smooth. It is free, and
$Z/G_0=N$. Thus, $Z\to N$ (\ref{g105}) is a principal bundle with
the structure group $G_0=G$ (\ref{g100}).

Furthermore, this principal bundle over a paracompact smooth
manifold $N$ is trivial as follows. In accordance with the
well-known theorem,$^{32}$ its structure group $G$ (\ref{g100}) is
reducible to the maximal compact subgroup $T^r$, which is also the
maximal compact subgroup of the group product
$\op\times^rGL(1,\Bbb C)$. Therefore, the equivalence classes of
$T^r$-principal bundles $\xi$ are defined as
\be
c(\xi)=c(\xi_1\oplus\cdots\oplus \xi_r)=(1+c_1(\xi_1))\cdots
(1+c_1(\xi_r))
\ee
by the Chern classes $c_1(\xi_i)\in H^2(N,\Bbb Z)$ of
$U(1)$-principal bundles $\xi_i$ over $N$.$^{32}$ Since the
cohomology group $H^2(N,\Bbb Z)$ of $N$ is trivial, all Chern
classes $c_1$ are trivial, and the principal bundle $Z\to N$ is
also trivial. This principal bundle can be provided with the
following coordinate atlas.

Let us consider the fibered manifold $S:Z\to W$ (\ref{g106}).
Because functions $S_\la$ are constant on fibers of the fiber
bundle $Z\to N$ (\ref{g105}), the fibered manifold (\ref{g106})
factorizes through the fiber bundle (\ref{g105}), and we have the
composite fibered manifold (\ref{g107}). Let us provide the
principal bundle $Z\to N$ with a trivialization
\mar{g110}\beq
Z=N\times \Bbb R^{m-r}\times T^r\to N, \label{g110}
\eeq
whose fibers are endowed with the standard coordinates
$(y^\la)=(t^a,\vf^i)$ on the toroidal cylinder (\ref{g100}). Then
the composite fibered manifold (\ref{g107}) is provided with the
fibered coordinates
\mar{g108}\ben
&& (J_\la,x^A,t^a,\vf^i), \label{g108}\\
&& \la=1,\ldots, m, \quad A=1, \ldots, 2(n-m), \quad a=1, \ldots,
m-r, \quad i=1,\ldots, r, \nonumber
\een
where $J_\la =S_\la(z)$ are coordinates on the base $W$ induced by
Cartesian coordinates on $\Bbb R^m$, and $(J_\la, x^A)$ are
fibered coordinates on the fibered manifold $\zeta:N\to W$. The
coordinates $J_\la$ on $W\subset \Bbb R^m$ and the coordinates
$(t^a,\vf^i)$ on the trivial bundle (\ref{g110}) possess the
identity transition functions, while the transition function of
coordinates $(x^A)$ depends on the coordinates $(J_\la)$ in
general.

The Hamiltonian vector fields $v_\la$ on $Z$ relative to the
coordinates (\ref{g108}) take the form
\mar{ww25}\beq
v_\la=v_\la^a(x)\dr_a + v^i_\la(x)\dr_i. \label{ww25}
\eeq
Since these vector fields commute (i.e., fibers of $Z\to N$ are
isotropic), the symplectic form $\Om$ on $Z$ reads
\mar{g103}\beq
\Om=\Om^\al_\bt dJ_\al\w dy^\bt + \Om_{\al A}dy^\al\w dx^A +
\Om^{\al\bt} dJ_\al\w dJ_\bt + \Om^\al_A d J_\al\w dx^A +\Om_{AB}
dx^A\w dx^B. \label{g103}
\eeq

\begin{lem} \label{g144} \mar{g144}
The symplectic form $\Om$ (\ref{g103}) is exact.
\end{lem}

\begin{proof} In accordance with the well-known K\"unneth formula,
the de Rham cohomology group of the product (\ref{g110}) reads
\be
H^2(Z)=H^2(N)\oplus H^1(N)\ot H^1(T^r) \oplus H^2(T^r).
\ee
By the de Rham theorem,$^{32}$ the de Rham cohomology $H^2(N)$ is
isomorphic to the cohomology $H^2(N,\Bbb R)$ of $N$ with
coefficients in the constant sheaf $\Bbb R$. It is trivial since
$H^2(N,\Bbb R)=H^2(N,\Bbb Z)\ot\Bbb R$ where $H^2(N,\Bbb Z)$ is
trivial. The first cohomology group $H^1(N)$ of $N$ is trivial
because $N$ is simply connected. Consequently, $H^2(Z)=H^2(T^r)$.
Then the closed form $\Om$ (\ref{g103}) is exact since it does not
contain the term $\Om_{ij}d\vf^i\w d\vf^j$.
\end{proof}

Thus, we can write
\mar{g113}\beq
\Om=d\Xi, \qquad \Xi=\Xi^\la(J_\al,x^B,y^\al) dJ_\la +
\Xi_\la(J_\al,x^B) dy^\la +\Xi_A(J_\al,x^B,y^\al) dx^A.
\label{g113}
\eeq
Up to an exact summand, the Liouville form $\Xi$ (\ref{g113}) is
brought into the form
\mar{g114}\beq
\Xi=\Xi^\la(J_\al,x^B,y^\al) dJ_\la + \Xi_i(J_\al,x^B) d\vf^i
+\Xi_A(J_\al,x^B,y^\al) dx^A, \label{g114}
\eeq
i.e., it does not contain the term $\Xi_a dt^a$.

The Hamiltonian vector fields $v_\la$ (\ref{ww25}) obey the
relations $v_\la\rfloor\Om=-dJ_\la$, which result in the
coordinate conditions
\mar{ww22}\beq
  \Om^\al_\bt \vt^\bt_\la=\dl^\al_\la, \qquad \Om_{A\bt}\vt^\bt_\la=0.
\label{ww22}
\eeq
The first of them shows that $\Om^\al_\bt$ is a nondegenerate
matrix independent of coordinates $y^\la$. Then the second one
implies $\Om_{A\bt}=0$.

Since $\Xi_a=0$ and $\Xi_i$ are independent of $\vf^i$, it follows
from the relations
\be
\Om_{A\bt}=\dr_A\Xi_\bt-\dr_\bt\Xi_A=0
\ee
that $\Xi_A$ are independent of coordinates $t^a$ and at most
affine in $\vf^i$. Since $\vf^i$ are cyclic coordinates,
  $\Xi_A$ are independent of $\vf^i$. Hence,
$\Xi_i$ are independent of coordinates $x^A$, and the Liouville
form $\Xi$ (\ref{g114}) reads
\mar{ac2}\beq
\Xi=\Xi^\la(J_\al,x^B,y^\al)dJ_\la + \Xi_i(J_\al) d\vf^i
+\Xi_A(J_\al,x^B)dx^A. \label{ac2}
\eeq
Because entries $\Om^\al_\bt$ of $d\Xi=\Om$ are independent of
$y^\la$, we obtain the following.

(i) $\Om^\la_i=\dr^\la\Xi_i-\dr_i\Xi^\la$. Consequently,
$\dr_i\Xi^\la$ are independent of $\vf^i$, and so are $\Xi^\la$
since $\vf^i$ are cyclic coordinates. Hence,
$\Om^\la_i=\dr^\la\Xi_i$ and $\dr_i\rfloor\Om=-d\Xi_i$. A glance
at the last equality shows that $\dr_i$ are Hamiltonian vector
fields. It follows that, from the beginning, one can separate $r$
integrals of motion, say $H_i$ again, whose Hamiltonian vector
fields are tangent to invariant tori. In this case, the
Hamiltonian vector fields $v_\la$ (\ref{ww25}) read
\mar{ww25'}\beq
v_a=\dr_a, \qquad v_i=v_i^k(x)\dr_k, \label{ww25'}
\eeq
where the matrix function $v_i^k(x)$ is nondegenerate. Moreover,
the coordinates $t^a$ are exactly the flow parameters $s^a$.
Substituting the expressions (\ref{ww25'}) into the first
condition (\ref{ww22}), we obtain
\be
\Om=dJ_a\w ds^a + (v^{-1})^i_k dJ_i\w d\vf^k +
\Om^{\al\bt}dJ_\al\w dJ_\bt + \Om_A^\la dJ_\la\w dx^A
+\Om_{AB}dx^A\w dx^B.
\ee
It follows that $\Xi_i$ are independent of $J_a$, and so are
$(v^{-1})^k_i=\dr^k\Xi_i$.

(ii) $\Om^\la_a=-\dr_a\Xi^\la=\dl^\la_a$. Hence,
$\Xi^a=-t^a+E^a(J_\la,x^B)$ and $\Xi^i$ are independent of $t^a$.

In view of items (i) -- (ii), the Liouville form $\Xi$ (\ref{ac2})
reads
\be
\Xi=(-t^a+E^a(J_\la,x^B))dJ_a + E^i(J_\la,x^B) dJ_i + \Xi_i(J_j)
d\vf^i + \Xi_A(J_\la,x^B)dx^A.
\ee
Since the matrix $\dr^k\Xi_i$ is nondegenerate, we can perform the
coordinate transformations $I_a=J_a$, $I_i=\Xi_i(J_j)$ together
with the coordinate transformations
\mar{g142}\beq
 t'^a = -t^a+E^a(J_\la,x^B),
\qquad \vf'^i =\vf^i-E^j(J_\la,x^B)\frac{\dr J_j}{\dr I_i}.
\label{g142}
\eeq
These transformations bring $\Om$ into the form (\ref{nc3'}).

\bigskip
\bigskip

\noindent {\bf III. PROOF OF THEOREM 5}
\bigskip

Given a fibration $H$ (\ref{nc4}), let $V$ be an open subset of
its base $N$ which satisfies condition (iii) of Theorem \ref{nc0},
i.e., there is a set $\{C_1,\ldots, C_m\}$ of $m$ independent
Casimir functions of the coinduced Poisson structure $\{,\}_N$ on
$V$. Note that such functions always exist around any point of
$N$. Let $Z_V$ be the restriction of the fiber bundle $Z\to N$
onto $V\subset Z$. By virtue of Lemma \ref{nc8}, $Z_V\to V$ is a
fibration in invariant submanifolds of a PIS $\{H^*C_\la\}$, where
$H^*C_\la$ are the pull-back of the Casimir functions $C_\la$ onto
$Z_V$.

Let $v_\la$ be Hamiltonian vector fields of functions $H^*C_\la$.
Since
\mar{j21}\beq
H^*C_\la(z)= (C_\la\circ H)(z)= C_\la(H_i(z)),\qquad z\in Z_V,
\label{j21}
\eeq
the Hamiltonian vector fields $v_\la$ restricted to any fiber $M$
of $Z_V$ are linear combinations of the Hamiltonian vector fields
$\vt_i$ of integrals of motion $H_i$. It follows that $v_\la$ are
elements of a finite-dimensional real Lie algebra of vector fields
on $M$ generated by the vector fields $\vt_i$. Since vector fields
$\vt_i$ are complete, the vector fields $v_\la$ on $M$ are also
complete.$^{28}$ Consequently, the Hamiltonian vector fields
$v_\la$ are complete on $Z_V$. Then the conditions of Theorem
\ref{nc6} for a PIS $\{H^*C_\la\}$ on the symplectic manifold
$(Z_V,\Om)$ hold.

In accordance with Theorem \ref{nc6}, we have a composite fibered
manifold
\mar{g150}\beq
Z_V\ar^H V\ar^C W, \label{g150}
\eeq
where $C:V\to W$ is a fibered manifold of level surfaces of the
Casimir functions $C_\la$. The fibered manifold (\ref{g150}) is
provided with the adapted fibered coordinates $(J_\la, x^A,
y^\la)$ (\ref{g108}), where $J_\la$ are values of the Casimir
functions and $(y^\la)=(t^a,\vf^i)$ are coordinates on a toroidal
cylinder. Since $C_\la=J_\la$ are Casimir functions on $V$, the
symplectic form $\Om$ (\ref{g103}) on $Z_V$ reads
\mar{g141}\beq
\Om=\Om^\al_\bt dJ_\al\w dy^\bt + \Om_{\al A}dy^\al\w dx^A +
\Om_{AB} dx^A\w dx^B. \label{g141}
\eeq
In particular, it follows that transition functions of coordinates
$x^A$ on $V$ are independent of coordinates $J_\la$, i.e., $C:V\to
W$ is a trivial bundle.

By virtue of Lemma \ref{g144}, the symplectic form (\ref{g141}) is
exact, i.e., $\Om=d\Xi$, where the Liouville form $\Xi$
(\ref{ac2}) is
\be
\Xi=\Xi^\la(J_\al,y^\la)dJ_\la + \Xi_i(J_\al) d\vf^i
+\Xi_A(x^B)dx^A.
\ee
It is brought into the form
\be
\Xi=(-t^a+E^a(J_\la))dJ_a + E^i(J_\la) dJ_i + \Xi_i(J_j) d\vf^i +
\Xi_A(x^B)dx^A.
\ee
Then the coordinate transformations (\ref{g142})
\mar{g151}\ben
&& I_a=J_a, \qquad I_i=\Xi_i(J_j), \nonumber \\
&& t'^a = -t^a+E^a(J_\la), \qquad \vf'^i
=\vf^i-E^j(J_\la)\frac{\dr J_j}{\dr I_i}, \label{g151}
\een
bring $\Om$ (\ref{g141}) into the form (\ref{nc3}). In comparison
with the general case (\ref{g142}), the coordinate transformations
(\ref{g151}) are independent of coordinates $x^A$. Therefore, the
angle coordinates $\vf'^i$ possess identity transition functions
on $V$.

Theorem \ref{nc0} restarts Theorem \ref{nc0'} if one considers an
open subset $U$ of $V$ admitting the Darboux coordinates $x^A$ on
the symplectic leaves of $U$.

The proof of Theorem \ref{nc0} gives something more. Let $\cH$ be
a Hamiltonian of a CIS. It is the pull-back onto $Z_V$ of some
Casimir function on $V$. Since $(I_\la,x^A)$ are coordinates on
$V$, they are also integrals of motion of $\cH$. Though the
original integrals of motion $H_i$ are smooth functions of
coordinates $(I_\la,x^A)$, the Casimir functions (\ref{j21})
\be
C_\la(H_i(I_\m,x^A))= C_\la(I_\m)
\ee
and, in particular, a Hamiltonian $\cH$ depend only on the action
coordinates $I_\la$. Hence, the equations of motion of a CIS take
the form
\be
\dot y'^\la=\frac{\dr\cH}{\dr I_\la}, \quad I_\la={\rm const.},
\quad x^A={\rm const.}
\ee

\end{document}